\documentclass[
    ,final            
  ]
  {aipproc}

\layoutstyle{8x11single}
\usepackage{graphicx}
\usepackage{rotating}

\def\RR{\mbox{$\rm I\!R$}}
\def\pmatrix{ \left( \begin{array} }
\def\endpmatrix{ \end{array} \right) }
\def\P{{\mathcal P}}
\def\O{{\Omega}}

\begin{document}

\title{Numerical comparisons among some methods for Hamiltonian problems}


\author{Luigi Brugnano}{
  address={Dipartimento di Matematica ``U.\,Dini'', Viale Morgagni 67/A, 50134
Firenze, Italy.}
}

\author{Felice Iavernaro}{
  address={Dipartimento di Matematica, Via Orabona 4, 70125 Bari, Italy.}
}

\author{Donato Trigiante}{
  address={Dipartimento di Energetica ``S.\,Stecco'', Universit\`a di
Firenze, Via Lombroso 6/17, 50134 Firenze (Italy).} }


\maketitle



We here report a few numerical tests comparing geometric
integrators, of Runge-Kutta type, described by Butcher tableaus in
the following form: \begin{equation}\label{bt} {\rm
(a)}\quad\begin{array}{c|c}
\begin{array}{c} c_1\\ \vdots\\ c_s\end{array} & \P X_s(\alpha) \P^{-1}
\\ \hline & b_1~\dots~ b_s
\end{array}\qquad \qquad \qquad {\rm (b)}\quad \begin{array}{c|c}
\begin{array}{c} c_1\\ \vdots\\ c_k\end{array} & \P_k \hat{X}_s
\P_k^T\O
\\ \hline & b_1~\dots~ b_k
\end{array}
\end{equation}
where $k\ge s$, ~$\{c_1<c_2<\dots<c_\ell\}$~ and ~$\{b_1,\dots,
b_\ell\}$~ are the abscissae and the weights of the Gauss-Legendre
quadrature formula in the interval $[0,1]$, $\ell=s,k$,
$$
X_s(\alpha) = \pmatrix{cccc}
\frac{1}2 & -\xi_1 &&\\
\xi_1     &0      &\ddots&\\
          &\ddots &\ddots    &-(\xi_{s-1}+\alpha)\\
          &       &\xi_{s-1}+\alpha &0\\
\endpmatrix, \quad
\hat{X}_s = \pmatrix{ccccc}
\frac{1}2 & -\xi_1 &&&\\
\xi_1     &0      &\ddots&&\\
          &\ddots &\ddots    &-\xi_{s-1}&\\
          &       &\xi_{s-1} &0 &\\
          &       &          &\xi_s&0
\endpmatrix,
\quad \xi_j=\frac{1}{2\sqrt{4j^2-1}},$$
$\O={\rm diag}(b_1,\dots,b_k)$ and, finally, by considering the
Legendre polynomials $P_j(\tau)$ of degree $j-1$, for $j\ge1$,
shifted and normalized in the interval $[0,1]$ so that $\int_0^1
P_i(\tau)P_j(\tau) \mathrm{d} \tau = \delta_{ij}$ (the Kronecker
symbol), $\P=(P_j(c_i))\in\RR^{s\times s}$,
$\P_k=(P_j(c_i))\in\RR^{k\times s+1}$. Method (\ref{bt})-(a)
reduces to the $s$-stage Gauss-Legendre method when $\alpha=0$
(see, e.g., \cite[pp.\,77 ff.]{HW}). The same happens to method
(\ref{bt})-(b) when $k=s$ \cite{BIT2,BIT3}. The $s$-stage
Gauss-Legendre method is known to be a symplectic integrator of
order $2s$, able to preserve quadratic invariants for Hamiltonian
problems in canonical form \cite{HLW}. On the other hand, under
suitable mild assumptions \cite{BIT4} the parameter $\alpha$ in
(\ref{bt})-(a) can be tuned, at each step, in order to obtain also
the conservation of the Hamiltonian (see also \cite{BIT5}): let us
denote such methods by EQUIP$(s)$ ({\em E}nergy and {\em
QU}adratic {\em I}nvariants {\em P}reserving) methods. Finally,
the formulae (\ref{bt})-(b) define the class of HBVM$(k,s)$
methods \cite{BIT00,BIT0,BIT1,BIT2,BIT3}, able to preserve
polynomial Hamiltonian functions of degree $\nu$, provided that
$k\ge (\nu s)/2$ (obviously, a {\em practical} conservation of
energy is obtained, for all suitably regular Hamiltonian
functions, provided that $k$ is large enough). The order of all
the above mentioned methods is $2s$. In the following we fix
$s=3$.

\begin{figure}\caption{\label{fig1}}
\centerline{\includegraphics[width=12cm,height=6cm]{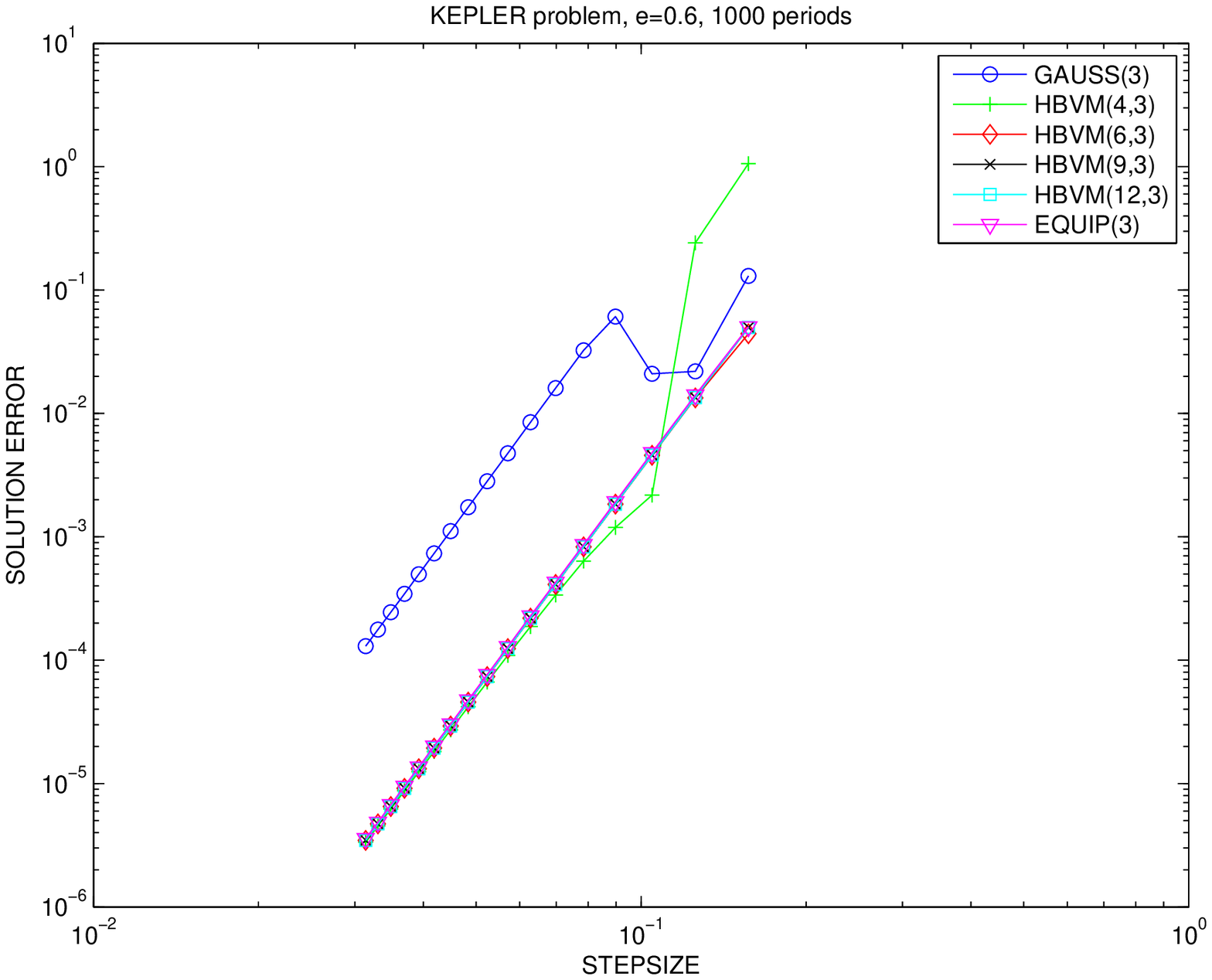}}
\end{figure}

\begin{figure}\caption{\label{fig2}}
\centerline{\includegraphics[width=12cm,height=6cm]{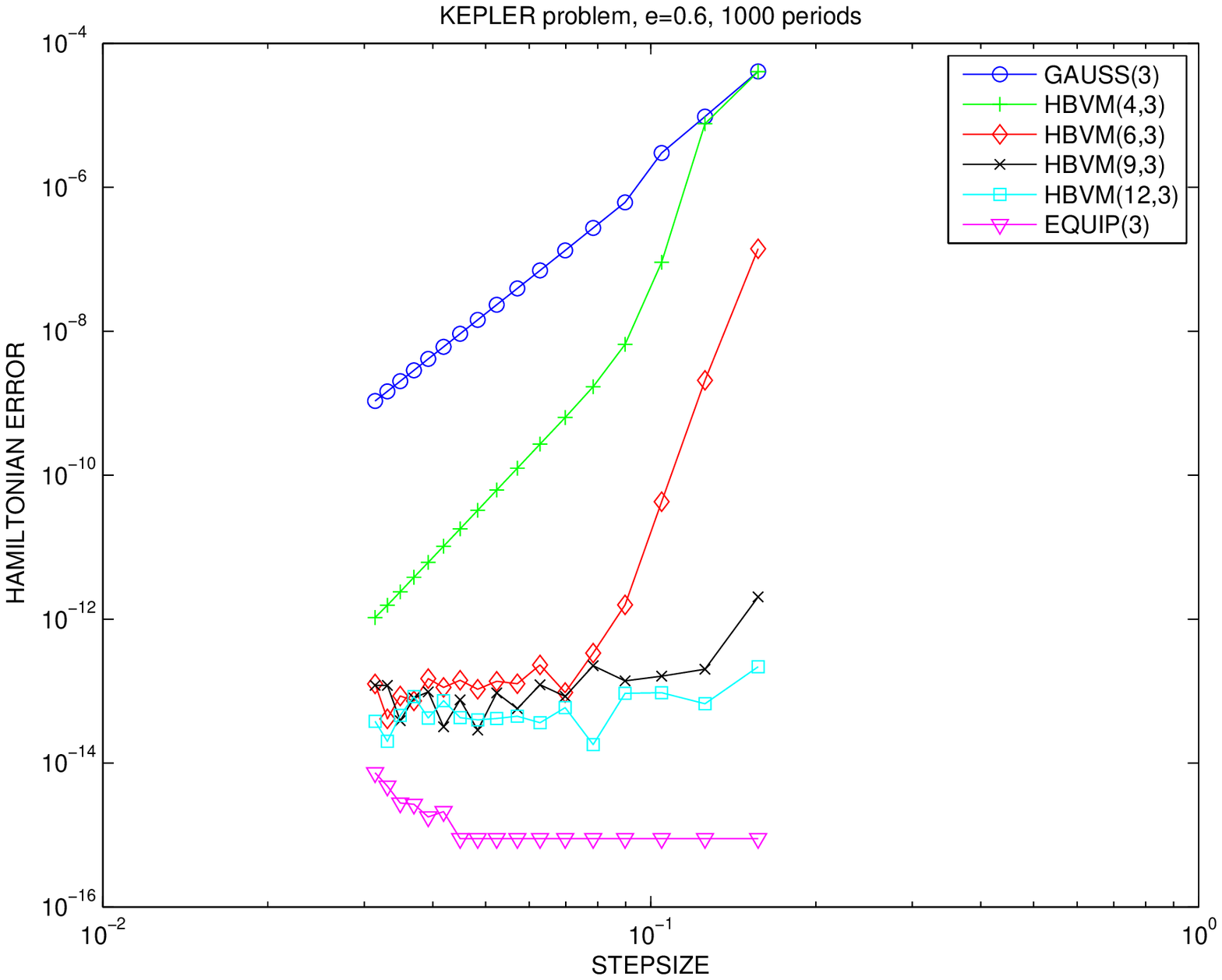}}
\end{figure}

\begin{figure}\caption{\label{fig3}}
\centerline{\includegraphics[width=12cm,height=6cm]{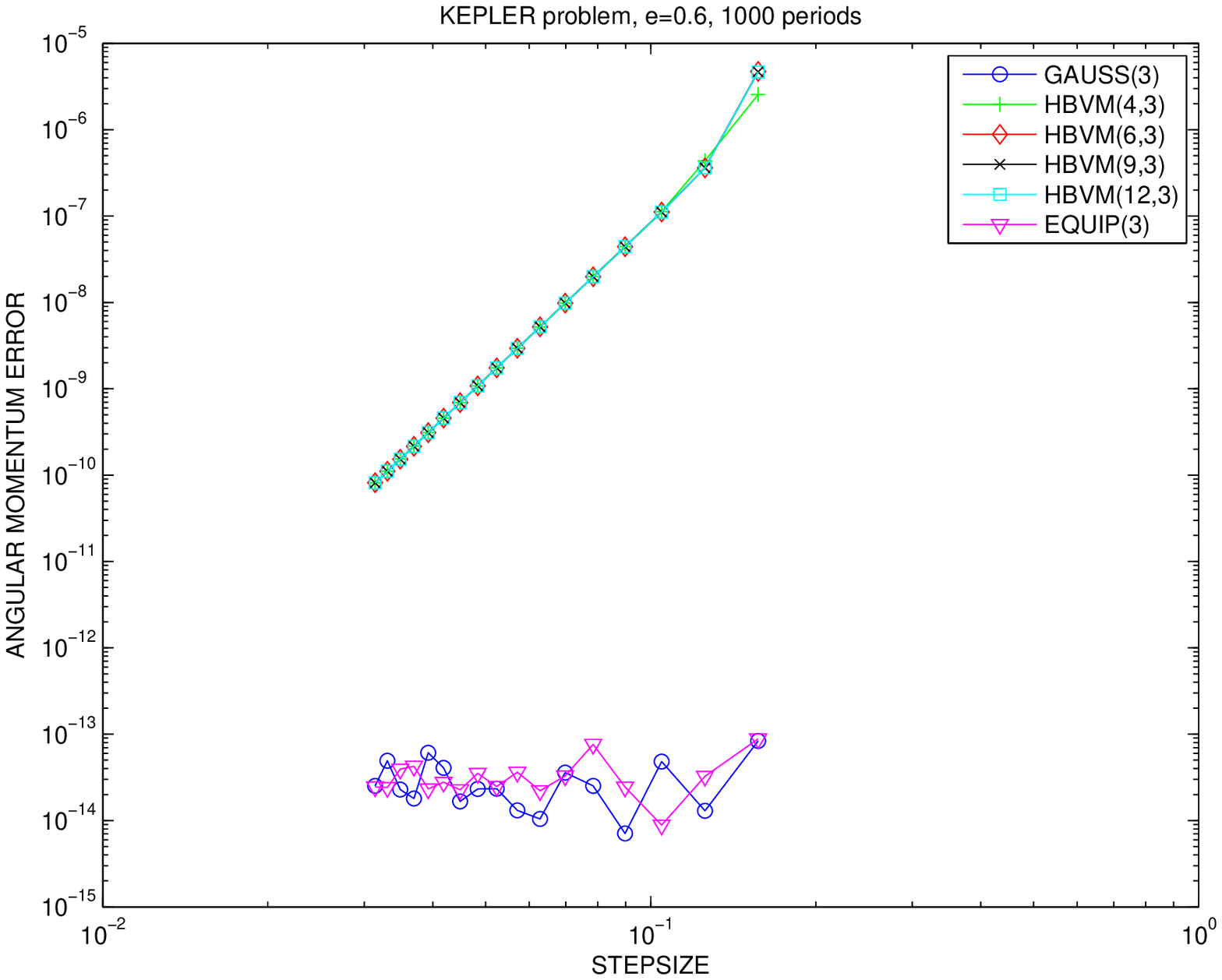}}
\end{figure}

In Figures~\ref{fig1}--\ref{fig3} we plot the errors (in the
solution, in the Hamiltonian, and in the angular momentum,
respectively) versus the (constant) stepsize used, for the
GAUSS(3) ($\equiv$ HBVM(3,3)), HBVM(4,3), HBVM(6,3), HBVM(9,3),
HBVM(12,3), and EQUIP(3) methods applied to the Kepler problem
\cite[pp.\,7--9]{HLW}, with eccentricity $e=0.6$, over 1000
periods. As one can see (Figure~\ref{fig1}), the order of all
methods is confirmed to be 6, even though the error constants of
HBV($k$,3), $k>3$, and EQUIP(3) methods turn out to be apparently
the same, and approximately 40 times less than that of GAUSS(3)
($\equiv$ HBVM(3,3)). The error in the Hamiltonian
(Figure~\ref{fig2}), as expected, decreases for HBVM($k$,3)
methods, as $k$ is increased (with order $2k$ \cite{BIT1}, until
round-off errors prevail), and practical conservation is obtained
for $k\ge9$. EQUIP(3) clearly conserves, by its own definition,
the Hamiltonian. Finally (Figure~\ref{fig3}), the error in the
angular momentum (which is a quadratic invariant) is negligible
for GAUSS(3) and EQUIP(3) methods, and decreasing at the same rate
$6$ $(\equiv2s)$ with the stepsize, for HBVM($k$,3), $k>3$,
methods. This is to be expected, since this error only depends on
matrix $\hat{X}_3$ (see (\ref{bt})-(b)), which is the same for all
such methods.

\begin{figure}\caption{\label{fig4}}
\centerline{\includegraphics[width=12cm,height=6cm]{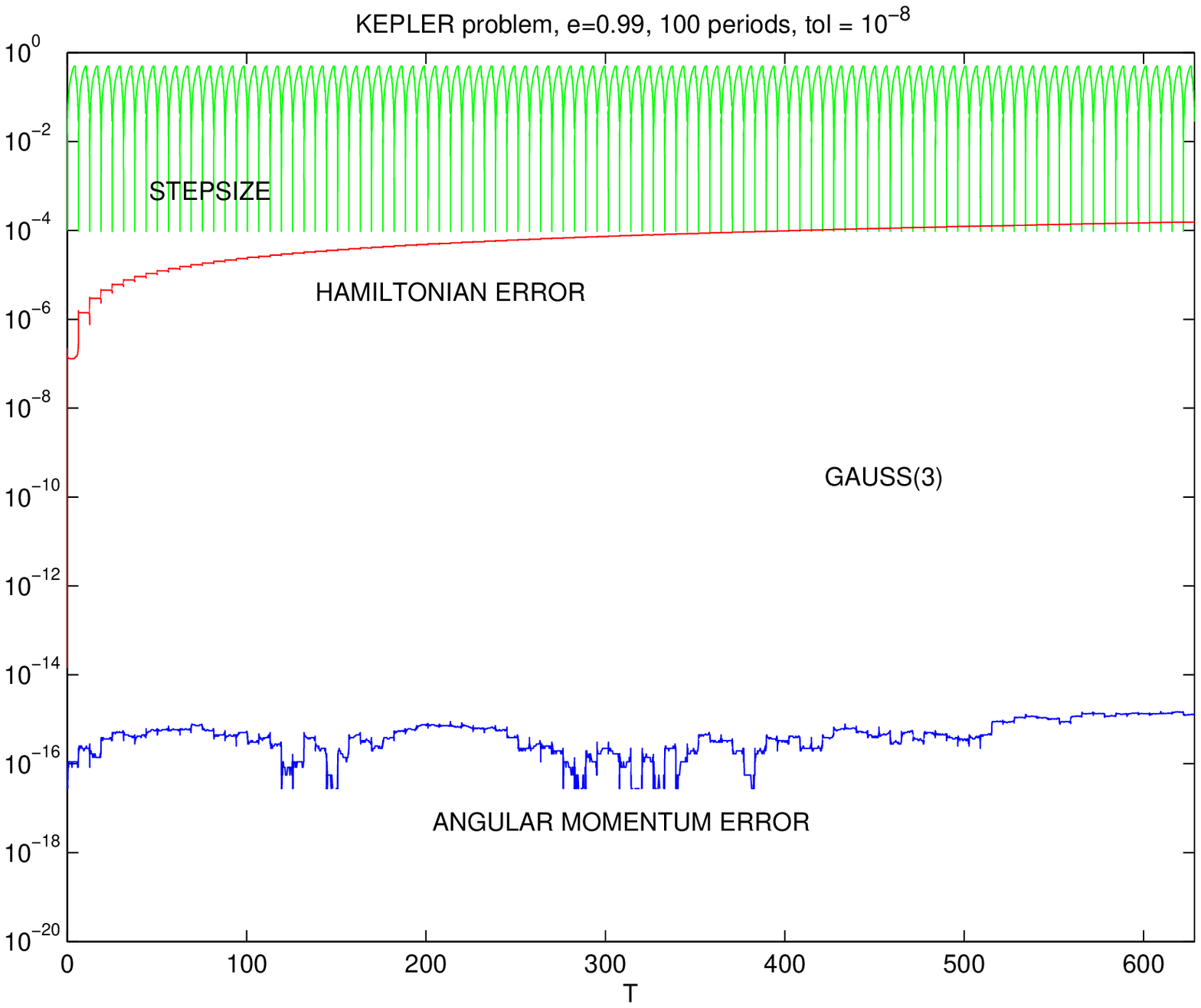}}
\end{figure}

\begin{figure}\caption{\label{fig5}}
\centerline{\includegraphics[width=12cm,height=6cm]{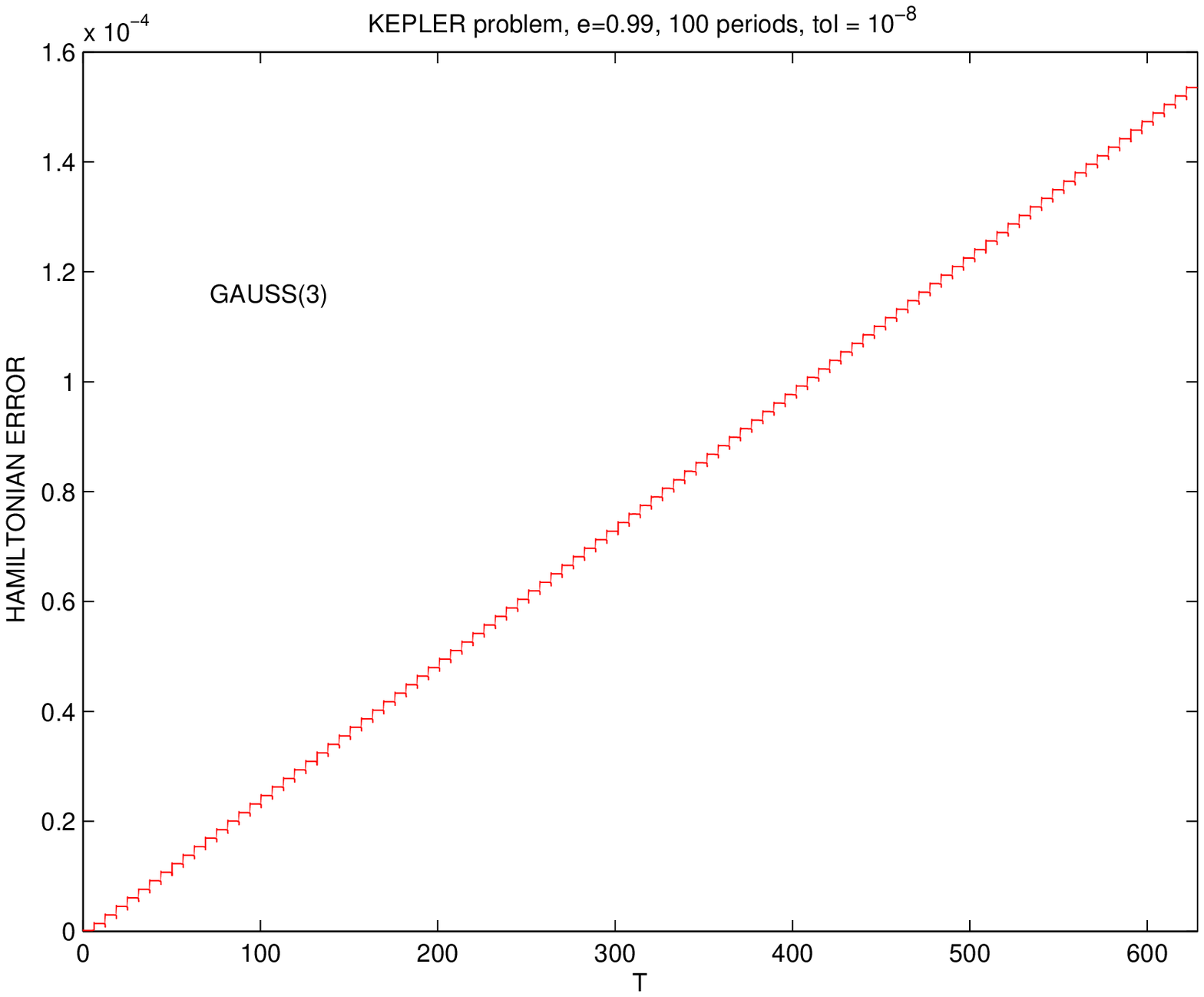}}
\end{figure}

\begin{figure}\caption{\label{fig6}}
\centerline{\includegraphics[width=12cm,height=6cm]{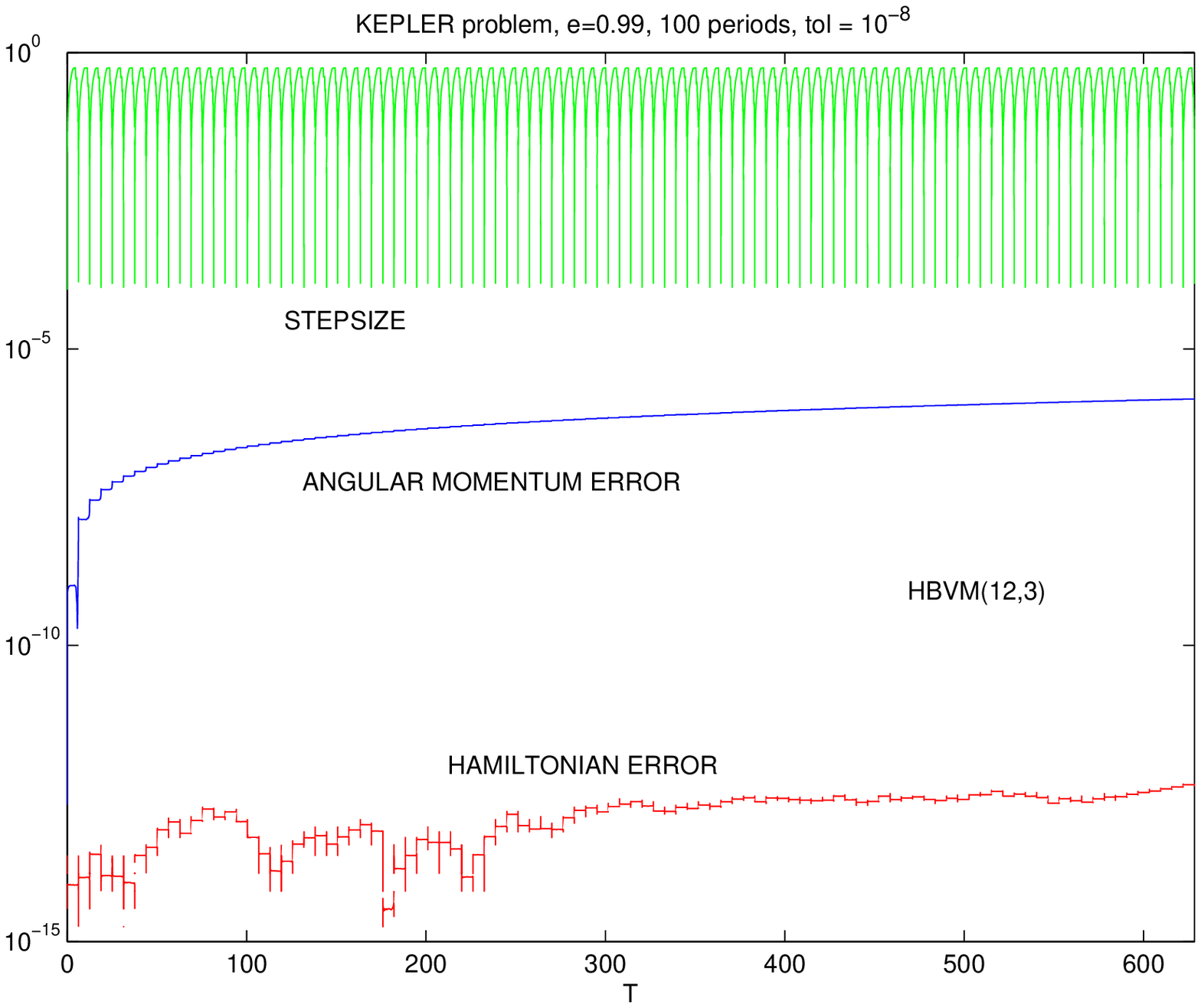}}
\end{figure}

\begin{figure}\caption{\label{fig7}}
\centerline{\includegraphics[width=12cm,height=6cm]{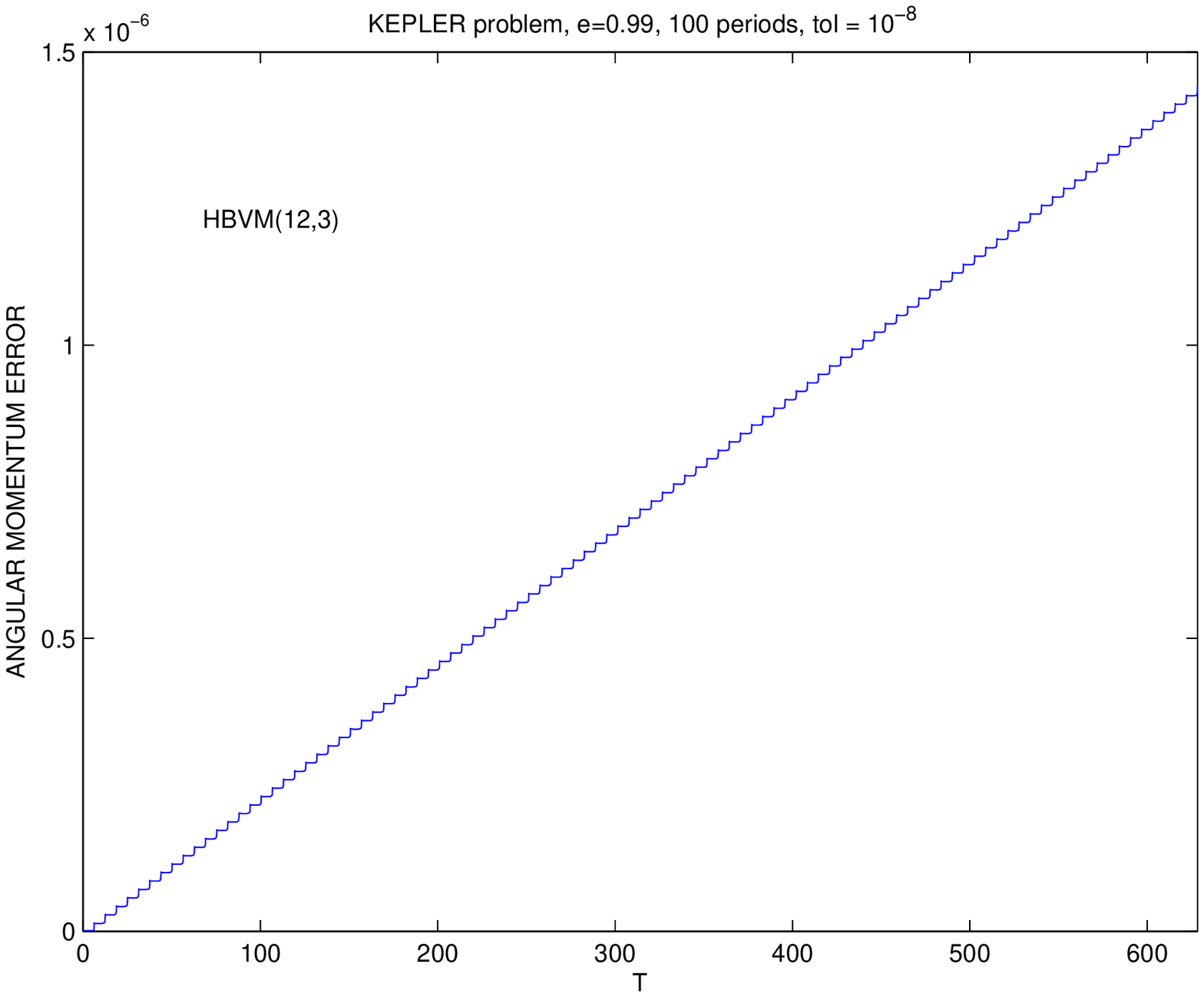}}
\end{figure}

\begin{figure}\caption{\label{fig8}}
\centerline{\includegraphics[width=12cm,height=6cm]{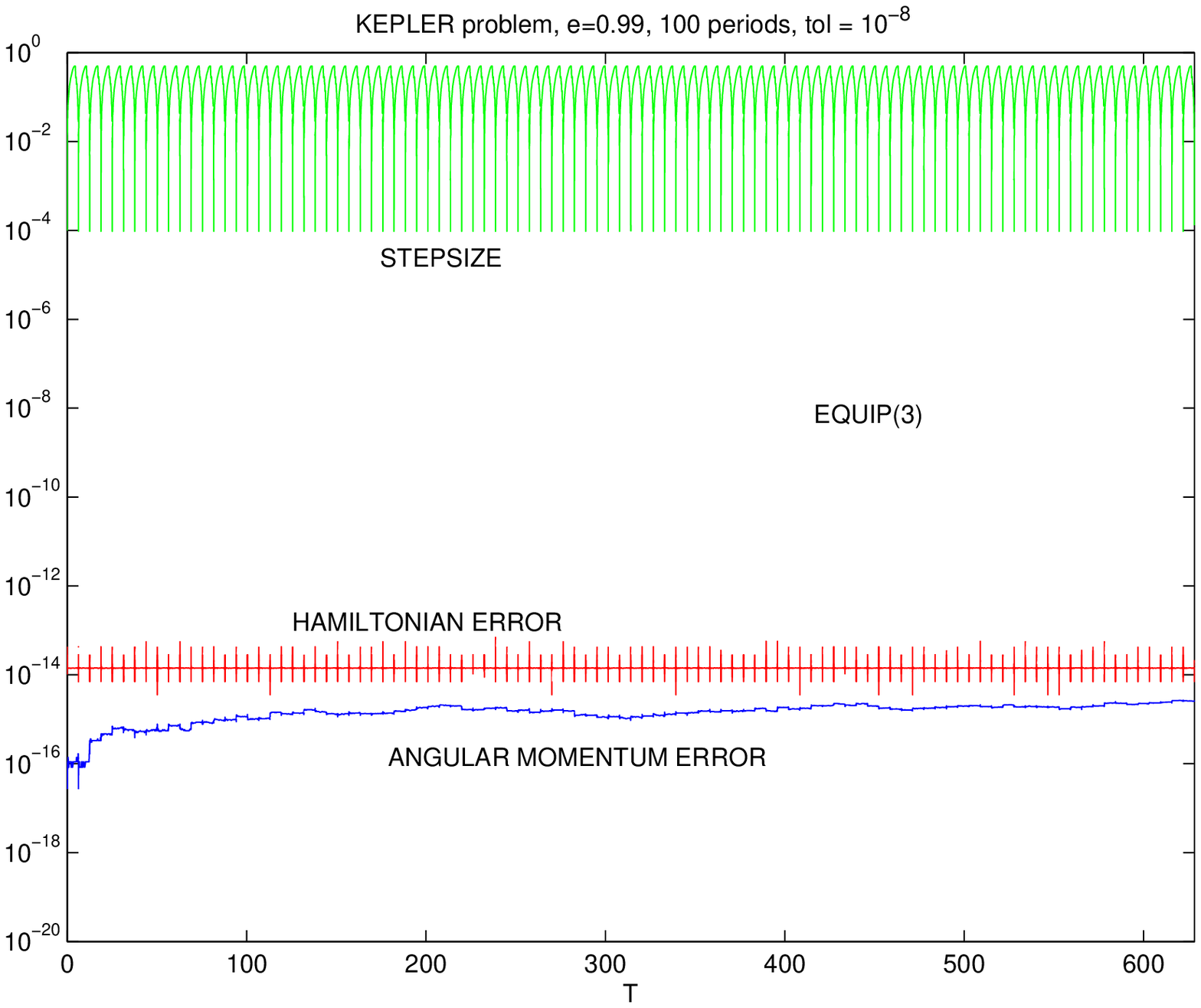}}
\end{figure}

To conclude, we report the numerical results, by using variable
stepsize with a \underline{standard} stepsize selection strategy
($tol=10^{-8}$) , for the GAUSS(3), HBVM(12,3), and EQUIP(3)
methods applied to the Kepler problem, with eccentricity $e=0.99$,
over 100 periods. All methods select stepsizes in the range
$10^{-4}\div 10^0$. As is well known \cite{HLW} standard stepsize
strategies don't work well with symplectic methods, so that
GAUSS(3), though preserving the angular momentum, exhibits a drift
in the numerical Hamiltonian (see Figures~\ref{fig4} and
\ref{fig5}). On the contrary, HBVM(12,3) practically conserves the
Hamiltonian but exhibits a drift in the angular momentum (see
Figures~\ref{fig6} and \ref{fig7}). At last, from
Figure~\ref{fig8} we conclude that only EQUIP(3) preserves
\underline{both} the energy and the angular momentum, when a
standard mesh selection strategy is used.


\bibliographystyle{aipproc}   



\end{document}